\documentclass[12pt]{amsart}

\usepackage[arrow,dvips,matrix,arrow,ps,color,line,curve,frame]{xy}

\usepackage{amssymb}

\advance\textheight1cm

\let\oldlabel=\label
\def\prellabel{\marginparsep=1em
    \def\label##1{\oldlabel{##1}\ifmmode\else\ifinner\else
         \marginpar{{\footnotesize\ \\ \tt
                    ##1}}\fi\fi}}

\advance\headheight1.15pt

\let\:\colon
\let\epsilon\varepsilon
\let\xto\xrightarrow
\let\phi=\varphi
\let\theta=\vartheta

\let\Bbb=\mathbb

\def\opn#1#2{\def#1{\operatorname{#2}}}
\opn\gp{gp} \opn\Max{Max} \opn\Ker{Ker} \opn\Coker{Coker}
\opn\Ext{Ext} \opn\conv{conv} \opn\chara{char} \opn\n{n} \opn\h{h}
\opn\GL{GuL} \opn\SL{SL} \opn\sn{sn} \opn\inte{int} \opn\End{End}
\opn\rank{rank} \opn\Aff{Aff} \opn\Spec{Spec} \opn\Proj{Proj}
\opn\QF{QF} \opn\I{Im} \opn\Hom{Hom} \opn\Aut{Aut} \opn\W{Witt}
\opn\w{w} \opn\inte{int} \opn\pyr{pyr} \opn\l{l} \opn\r{r}
\opn\const{const}

\def\Pp{{\Bbb P}}
\def\ZZ{{\Bbb Z}}

\def\NN{{\Bbb N}}
\def\QQ{{\Bbb Q}}

\opn\End{End}
\opn\ch{ch}%

\def\Q{{\Box\kern1pt}}%

\def\k{{\bf k}}

\def\1{^{-1}}

\newtheorem{lemma}{Lemma}

\newtheorem{proposition}[lemma]{Proposition}

\theoremstyle{definition}

\begin{document}

\title[The coefficient field in the Nilpotence conjecture]{The coefficient
field in the Nilpotence\\ Conjecture for toric varieties}

\author{Joseph Gubeladze}

\thanks{Supported by MSRI, INTAS grant
99-00817 and TMR grant ERB FMRX CT-97-0107}


\address{A. Razmadze Mathematical Institute, Alexidze St. 1, 380093
Tbilisi, Georgia}

\address{Department of Mathematics, San Francisco
State University, San Francisco, CA 94132, USA}

\email{soso@math.sfsu.edu}

\begin{abstract}
The main result of the work ``The nilpotence conjecture in
$K$-theory of toric varieties'' is extended to \emph{all}
coefficient fields of characteristic 0, thus covering the class of
genuine toric varieties.
\end{abstract}

\maketitle

\section{The statement}\label{statement}

Let $R$ be a (commutative) regular ring, $M$ be arbitrary
commutative, cancellative, torsion free monoid without nontrivial
units, and $i$ be a nonnegative integral number. The {\em
nilpotence conjecture} asserts that for any sequence ${\bf
c}=(c_1,c_2,\dots)$ of natural numbers $\geq2$ and any element
$x\in K_i(R[M])$ there exists an index $j_x\in\NN$ such that
$(c_1\cdots c_j)_*(x)\in K_i(R)$ for all $j>j_x$.

Here $R[M]$ is the monoid $R$-algebra of $M$ and for a natural
number $c$ the endomorphism of $K_i(R[M])$, induced by the
$R$-algebra endomorphism $R[M]\to R[M]$, $m\mapsto m^c$, $m\in M$,
is denoted by $c_*$ (writing the monoid operation
multiplicatively). We call this action the \emph{multiplicative
action} of $\NN$ on $K_i(R[M])$.

Theorem 1.2 in \cite{G2} verifies the conjecture for the
coefficient rings $R$ of type $S^{-1}\k[T_1,\ldots,T_d]$ where
$\k$ is a number field and $S\subset\k[T_1,\ldots,T_d]$ is
arbitrary multiplicative subset of nonzero polynomials
($d\in\ZZ_+$). In particular, the nilpotence conjecture is valid
for purely transcendental extensions of $\QQ$. On the other hand
$K$-groups commute with filtered colimits. Therefore, the
following induction proposition, together with \cite[Theorem
2.1]{G2}, proves the conjecture for all characteristic 0 fields.

\begin{proposition}\label{finitext}
The validity of the nilpotence conjecture for a field
$\k_1$ of characteristic 0 transfers to any finite field extension
$\k_1\subset\k_2$.
\end{proposition}

This note should be viewed as an addendum to the main paper
\cite{G2}. In Sections \ref{ratgaldesc} and \ref{stienstra} for
the reader's convenience we recall in a convenient form the needed
results on Galois descent and Bloch-Stienstra operations in
$K$-theory.

\section{Galois descent for rational $K$-theory}\label{ratgaldesc}

We need the fact that rational $K$-theory of rings (essentially)
satisfies Galois descent -- a very special case of Thomason's
\'etale descent for localized versions of $K$-theory, first proved
in \cite{T} in the smooth case and then extended to the singular
case as an application of the new local-to-global technique
\cite{TT}. The argument below follows closely \cite[Lemma 2.13]{T}
and \cite[Proposition 11.10]{TT}. It is, of course, important that
Thomason's higher $K$-groups agree with those of Quillen in the
affine (or, more generally, quasiprojective) case \cite[Theorem
7.6 ]{TT}.

\begin{lemma}\label{galois}
Let $A\subset B$ be a finite Galois extension of noetherian rings
with the Galois group $G$. Assume $[B]=n[A]$ in $K_0(A)$ where,
naturally, $n=\#G$ (e.~g. B is a free $A$-module). Then
$K_i(A)\otimes\QQ=H^0(G,K_i(B)\otimes\QQ)$.
\end{lemma}

\begin{proof}
The push-out diagram of rings
\begin{equation}\label{1}
\xymatrix{B\ar[r]^{\iota_1}&B\otimes_AB\\
A\ar[r]_{\iota}\ar[u]^\iota&B\ar[u]_{\iota_2}}
\end{equation}
defines a pull-back diagram of the corresponding affine schemes in
the category of schemes. The latter diagram satisfies \emph{all}
the conditions of Proposition 3.18 in \cite{TT}. Therefore, by the
mentioned proposition we have the equality of the endomorphisms
\begin{equation}\label{2}
\iota_*\circ\iota^*=(\iota_2)^*\circ(\iota_1)_* :K_i(B)\to K_i(B)
\end{equation}
where $-_*$ refers to the functorial homomorphisms of the
$K_i$-groups and $-^*$ refers to the corresponding transfer maps
(contrary to the scheme-theoretical notation in
\cite{Q}\cite{T}\cite{TT}). Galois theory identifies the diagram
(\ref{1}) with the diagram
\begin{equation}\label{3}
\xymatrix{B\ar[r]^\Delta&B^n\\
A\ar[r]_{\iota}\ar[u]^{\iota}&B\ar[u]_{\Delta_G}}
\end{equation}
where $\Delta$ is the diagonal embedding and
$\bigl(\Delta_G(b)\bigr)_g=g(b)$. By the elementary properties of
$K$-groups \cite[\S2]{Q} and the equalities $(g^{-1})_*=g^*$,
$g\in G$ the equality (\ref{2}) and the diagram (\ref{3}) imply
$\iota_*\circ\iota^*=\sum_G(g^{-1})_*=\sum_Gg_*$. The other
composite $\iota^*\circ\iota_*$ is the multiplication by $n$ on
$K_i(A)$. Therefore, the homomorphisms $\iota_*\otimes\QQ$ and the
corresponding restriction of  $n^{-1}(\iota^*\otimes\QQ)$
establish the desired isomorphism.
\end{proof}

\section{Verschiebung and Frobenius}\label{stienstra}

We also need the Bloch-Stienstra formula on the relationship
between \emph{Verschiebung} and \emph{Frobenius} in the context of
the action of big Witt vectors on Nil-$K$-theory.

For a (commutative) ring $R$ the additive group of $\W(R)$ is the
multiplicative group $1+TR[[T]]$, $T$ a variable. The decreasing
filtration of subgroups $I_m(R)\subset 1+T^mR[[T]]$ makes $\W(R)$
a topological group and any element $\omega(T)\in\W(R)$ has a
unique convergent expansion $\omega(T)=\Pi_{n\geq1}(1-r_nT^n)$,
$r_n\in R$. When $\omega(T)\in I_m(R)$ then the expansion is of
the type $\omega(T)=\Pi_{n\geq m}(1-r_nT^n)$. The multiplicative
structure on $\W(R)$ is the unique continuous extension to the
whole $\W(R)$ of the pairing
$$
(1-rT^m)\star(1-sT^n)=(1-r^{n/d}s^{m/d}T^{mn/d})^d, \quad r,s\in
R,\quad d=\gcd(m,n).
$$
The assignment
$$
F_m:1-rT^n\mapsto(1-r^{n/d}T^{n/d})^d,\quad d=\gcd(m,n),\quad r\in
R
$$
extends to a (unique) ring endomorphism $F_m:\W(R)\to\W(R)$ which
is called the \emph{Frobenius endomorphism}.

When $\QQ\subset R$ the \emph{ ghost isomorphism} between the
multiplicative and additive groups
$$
-T\cdot\frac{d(\log)}{dT}:1+TR[[T]]\to
TR[[T]],\quad\alpha\mapsto\frac{-T}{\alpha}\cdot\frac{d\alpha}{dt}\\
$$
is actually a ring isomorphism  $\W(R)\to\Pi_1^\infty R$ where the
right hand side is viewed as $TR[[T]]$ under the assignment
$(r_1,r_2\ldots)\mapsto r_1T+r_2T^2+\cdots$. What we need here is
the fact that the diagonal injection $\lambda_T:R\to\Pi_1^\infty
R=\W(R)$ is invariant under $F_m$ for all $m\in\NN$.

The $R$-algebra endomorphism $R[T]\to R[T]$, $T\mapsto T^m$
induces a group endomorphism $V_m:NK_i(R)\to NK_i(R)$ -- the
\emph{Verschiebung} (see \cite[Theorem 4.7]{S} for such an
identification of $V_m$, in \cite{S} this map is first defined in
terms of the category $\text{\bf Nil}(R)$ of nilpotent
endomorphisms).

Bloch \cite{B} and then, in a systematic way, Stienstra \cite{S}
defined a $\W(R)$-module structure on $NK_i(R)$. The action of
$1-rT^m\in\W(R)$, $r\in R$ on $NK_i(R)$ is the effect of the
composite functor
$$
\Pp(R[T])\xto{t_m}\Pp(R[T])\xto{\bar r}\Pp(R[T])\xto{v_m}\Pp(R[T])
$$
where $v_m$ corresponds to the base change through $T\mapsto T^m$
(it gives rise to the $m$th Verschiebung), $\bar r$ corresponds to
the base change through $T\mapsto rT$, and $t_n$ corresponds to
the scalar restriction through $T\mapsto T^m$ (the corresponding
endomorphism of $NK_i(R)$ is the transfer). This determines the
action of the whole $\W(R)$ because any element of $NK_i(R)$ is
annihilated by the ideal $I_m(R)$ for some $m$ \cite[\S8]{S},
i.~e. $NK_i(R)$ is a \emph{continuous} $\W(R)$-module.

We have the following relationship between $F_m$ and $V_m$:
\begin{equation}\label{blst}
V_m(F_m(\alpha)\star z)=\alpha\star z,\quad\alpha\in\W(R),\quad
z\in NK_i(R),\quad m\in\NN,
\end{equation}
where $\star$ refers to the $\W(R)$-module structure. Actually,
(\ref{blst}) is proved in \cite[\S6]{S} for the elements $\alpha$
in the image of the Almkvist embedding $\ch:\tilde
K_0(\End(R))\to\W(R)$ \cite{A}. But the formula generalizes to the
whole $\W(R)$ due to the fact that $NK_i(R)$ is a continuous
$\W(R)$-module and the elements of type $1-rT^m$ are always in the
mentioned image.

In particular, when $\QQ\subset R$ the equality (\ref{blst})
implies

\begin{lemma}\label{frobenius}
For any natural number $m$ the Verschiebung $V_m$ is an $R$-linear
endomorphism of $NK_i(R)$.
\end{lemma}

Weibel \cite{W} has generalized the $\W(R)$-module structure to
the graded situation so that for a graded ring $A=A_0\oplus
A_1\oplus\cdots$  and a subring $R\subset A_0$ there is a {\it
functorial} continuous $\W(R)$-module structure on $K_i(A,A^+)$,
where $A^+=0\oplus A_1\oplus A_2\oplus\cdots$. (Although here we
restrict to the commutative case Weibel actually considers the
general situation when $R$ is in the center of $A$; it is exactly
the non-commutative case that is used in \cite{G2}). In
particular, the $\W(R)$-module structure on $NK_i(A,A^+)$ is the
restriction of that on $NK_i(A)$ under the embedding
$NK_i(A,A^+)\to NK_i(A)$ induced by the \emph{graded} homomorphism
\begin{align*}
\w:A_0\oplus A_1\oplus A_2\oplus\cdots\to A[T]=A+TA+T^2A+\cdots,\\
a_0\oplus a_1\oplus a_2\oplus\cdots\mapsto a_0+a_1T+a_2T^2+\cdots
\end{align*}
and, at the same time, the latter module structure is just the
scalar restriction of the $\W(A)$-module structure through
$\W(R)\to\W(A)$. That $NK_i(A,A^+)\to NK_i(A)$ is in fact an
embedding follows from the fact that $\w$ splits the non-graded
augmentation $A[T]\to A$, $T\mapsto 1$, in particular $K_i(A)\to
K_i(A[T])$ is a split monomorphism.

Below we will make use of the fact, due to the split exact
sequence $0\to A^+\to A\to A_0\to0$ ($A$ as above), that there is
a natural isomorphism $K_i(A,A^+)=K_i(A)/K_i(A_0)$. Also, the
latter group will be thought of as a direct summand of $K_i(A)$ in
a natural way.

\section{Proof of Proposition 1}\label{PROOF}

For clarity we let $c\mapsto c_*/\k$ denote the induced
multiplicative action of $\NN$ on $K_i(\k[M])/K_i(\k)$, $\k$ a
field.

Without loss of generality we can assume that $M$ is a finitely
generated monoid. Fix arbitrary embedding of $M$ into a free
monoid $\ZZ_+^r$ \cite[\S2.1]{BrG}. This gives rise to graded
structures on $\k_1[M]$ and $\k_2[M]$ so that the monoid elements
are homogeneous of positive degree. In view of the previous
section $K_i(\k_2[M])/K_i(\k_2)$ carries a $\k_2$-linear
structure, $K_i(\k_1[M])/K_i(\k_1)$ carries a $\k_1$-linear
structure, and the group homomorphism $K_i(\k_1[M])/K_i(\k_1)\to
K_i(\k_2[M])/K_i(\k_2)$ is $\k_1$-linear.

We can also assume that $\k_1\subset\k_2$ is a Galois extension.
In fact, if $\k_1\subset\k_3$ is a finite Galois extension such
that the conjecture is true for the monoid $\k_3$-algebras and
$\k_2\subset\k_3$ then the commutative squares
\begin{equation}\label{normal}
\xymatrix{K_i(\k_3[M])/K_i(\k_3)\ar[r]^{c_*/\k_3}&K_i(\k_3[M])/K_i(\k_3)\\
K_i(\k_2[M])/K_i(\k_2)\ar[r]_{c_*/\k_2}\ar[u]&K_i(\k_2[M])/K_i(\k_2)\ar[u]
},\qquad c\in\NN
\end{equation}
imply the validity of the nilpotence conjecture for $\k_2$ too
because the vertical homomorphisms in (\ref{normal}) are
monomorphisms. To see this observe that $\k_3[M]$ is a free module
over $\k_2[M]$ of rank $[\k_3:\k_2]$ and, hence, the composite of
the functorial map $K_i(\k_2[M])\to K_i(\k_3[M])$ with the
corresponding transfer map $K_i(\k_3[M])\to K_i(\k_2[M])$ is the
multiplication by $[\k_3:\k_2]$, which restricts to an
automorphism of the rational vector subspace
$K_i(\k_2[M])/K_i(\k_2)\subset K_i(\k_2[M])$.

Let $G$ be the Galois group of the extension $\k_1\subset\k_2$.
Then $\k_1[M]\subset\k_2[M]$ is a Galois extension of rings with
the same Galois group. The action of $G$ on $\W(\k_2)$ shows that
the $\k_1$-vector space $K_i(\k_2[M])/K_i(\k_2)$ is a Galois
$\k_1$-module. Then Lemma \ref{galois} implies
\begin{equation}\label{galdesc}
K_i(\k_2[M])/K_i(\k_2)=K_i(\k_1[M])/K_i(\k_1)\otimes\k_2.
\end{equation}
Since the nilpotence conjecture is valid for $\k_1$ every element
of $K_i(\k_1[M])/K_i(\k_1)$ is annihilated by high iterations of
the multiplicative action of $\NN$. It follows from
(\ref{galdesc}) and the commutative squares of type (\ref{normal})
for the extension $\k_1[M]\subset\k_2[M]$ that the same is true
for the elements of certain generating set of the $\k_2$-vector
space $K_i(\k_2[M])/K_i(\k_2)$. So we are done once it is shown
that $c_*/\k_2$ is a $\k_2$-linear map for every natural number
$c$.

Consider the commutative square of $\k_2$-algebras ($c\in\NN$)
\begin{equation}\label{linear}
\xymatrix{ \k_2[M][T]\ar[r]^{T\mapsto T^c}&\k_2[M][T]\\
\k_2[M]\ar[r]_{(-)^c}\ar[u]^{T^{\deg(-)}(-)^c}&\k_2[M]
\ar[u]_{T^{\deg(-)}(-)} }
\end{equation}
where every element $m\in M$ is fixed by the upper horizontal
homomorphism and
\begin{itemize}
\item
$\bigl((-)^c\bigr)(m)=m^c$,
\item
$\bigl(T^{\deg(-)}(-)^c\bigr)(m)=T^{\deg(m)}m^c$,
\item
$\bigl(T^{\deg(-)}(-)\bigr)(m)=T^{\deg(m)}m$.
\end{itemize}
Here $\deg(-)$ is the degree with respect to the fixed graded
structure on $\k_2[M]$. Thinking of $\k_2[M][T]$ as a graded ring
with respect to the powers of $T$ the vertical maps in
(\ref{linear}) become graded homomorphisms. We arrive at the
commutative diagram of groups
$$
\xymatrix{NK_i(\k_2[M])\ar[r]^{V_c}&NK_i(\k_2[M])\\
K_i(\k_2[M])/K_i(\k_2)\ar[r]_{c_*/\k_2}\ar[u]&K_i(\k_2[M])/K_i(\k_2)
\ar[u] }
$$
whose vertical maps are $\k_2$-linear by functoriality and the
upper horizontal map is $\k_2$-linear by Lemma \ref{frobenius}.
The left vertical homomorphism is actually a monomorphism because
the map $T^{\deg(-)}(-)$ in diagram (\ref{linear}) is of the type
$\w$, discussed at the end of Section \ref{stienstra}. These
conditions altogether imply the desired linearity. \qed

\end{document}